\documentclass[11pt]{llncs}
\usepackage {times}
\usepackage{amsfonts}
\usepackage{amssymb}
\usepackage{setspace}
\usepackage{lineno}
\newcommand {\bbox}{\rule{0.6em}{0.6em}}
\doublespace

\date{}

\title{Bounds on Isoperimetric Values of Trees}
\author {B. V. Subramanya Bharadwaj \thanks 
{Computer Science and Automation Department, Indian Institute
of Science, Bangalore- 560012 Email: subramanya@csa.iisc.ernet.in},  L. Sunil Chandran\thanks{Computer Science and Automation Department, Indian Institute of Science, Bangalore- 560012 Email: sunil@csa.iisc.ernet.in}} 
\institute{Indian Institute of Science, Bangalore, India}

\begin{document}
\maketitle
\pagestyle{plain}
\pagenumbering{arabic}
%\begin{pagewiselinenumbers}
\begin{abstract}
Let G = (V,E) be a finite, simple and undirected graph. For $S \subseteq V$, let $\delta(S,G) = \{ (u,v) \in E : u \in S \mbox { and } v \in V-S  \}$ be the edge boundary of $S$. Given an integer $i$, $1 \leq i \leq \vert V \vert$, let the edge isoperimetric value of $G$ at $i$ be defined as $b_e(i,G) = \min_{S \subseteq V; |S| = i} |\delta(S,G)|$. The edge isoperimetric peak of $G$ is defined as $b_e(G)=\max_{1 \leq j \leq \vert V \vert} b_e(j,G)$. Let $b_v(G)$ denote the vertex isoperimetric peak defined in a corresponding way.  The problem of determining a lower bound for the vertex isoperimetric peak in complete $t$-ary trees was recently considered in \cite{OatYam}. In this paper we provide bounds which improve those in \cite{OatYam}. Our results can be generalized to arbitrary (rooted) trees.
\paragraph{}
The depth $d$ of a tree is the number of nodes on the longest path starting from the root and ending at a leaf. In this paper we show that for a  complete binary tree of depth $d$ (denoted as $T_d^2$), $c_1d \leq b_e(T_d^2) \leq d$ and $c_2d \leq b_v(T_d^2) \leq d$ where $c_1$, $c_2$ are constants. For a complete $t$-ary tree of depth $d$ (denoted as $T_d^t$) and $d \geq c\log{t}$ where $c$ is a constant, we show that $c_1\sqrt{t}d \leq b_e(T_d^t) \leq td$ and $c_2\frac{d}{\sqrt{t}} \leq b_v(T_d^t) \leq d$ where $c_1$, $c_2$ are constants. At the heart of our proof we have the following theorem which works for an arbitrary rooted tree and not just for a complete $t$-ary tree. Let $T=(V,E,r)$ be a finite, connected and rooted tree - the root being the vertex $r$. Define a weight function $w:V \rightarrow \mathbb{N}$ where the weight $w(u)$ of a vertex $u$ is the number of its successors (including itself) and let the weight index $\eta(T)$ be defined as the number of distinct weights in the tree, i.e $\eta(T)=\vert \{ w(u):u \in V \} \vert$. For a positive integer $k$, let $\ell(k)=\vert \{ i \in \mathbb{N}:  1 \leq i \le \vert V \vert
\mbox {, }  b_e(i,G) \le k \} \vert$. We show that $\ell(k) \leq 2{{2\eta+k} \choose k}$.

\noindent {Keywords: } {\bf isoperimetric problem, binary trees, $t$-ary trees, pathwidth}.
 \end{abstract}
\section{Introduction}
Let $G=(V,E)$ be a simple, finite, undirected graph.
\begin{definition}
For $S \subseteq V$, the edge boundary  
$\delta(S,G)$ is  the set of edges of $G$ with exactly 
one end point in $S$. In other words,  
%\begin{linenomath*}
$$
\delta(S,G) = \{ (u,v) \in E : u \in S \mbox { and } v \in V-S  \} 
$$
%\end{linenomath*}
\end{definition}

\begin{definition}
For $S \subseteq V$, the vertex boundary $\phi(S,G)$ is defined similarly.
%\begin{linenomath*}
$$
\phi(S,G) = \{ v \in V -S: \exists u \in S,\mbox { such that } (u,v) \in E\} \nonumber
$$
%\end{linenomath*}
\end{definition}
\begin{definition}
Let $i$ be an integer where $1 \le i \le \vert V \vert$. For each $i$ define the edge isoperimetric value $b_e(i,G)$ and the vertex isoperimetric value $b_v(i,G)$ of $G$ at $i$ as follows
%\begin{linenomath*}
$$
b_e(i,G) = \min_{S \subseteq V; |S| = i} |\delta(S,G)|  \nonumber
$$
$$
b_v(i,G) = \min_{S \subseteq V; |S|=i} |\phi(S,G)| \nonumber
$$
%\end{linenomath*}
\end{definition}
\begin{definition}
For any graph $G$ define the edge and the vertex isoperimetric peaks $b_e(G)$, $b_v(G)$ as,
%\begin{linenomath*}
$$
b_e(G) = \max_{1 \leq i \leq \vert V \vert} b_e(i,G) \nonumber
$$
$$
b_v(G) = \max_{1 \leq i \leq \vert V \vert} b_v(i,G) \nonumber
$$
%\end{linenomath*}
\end{definition}

The edge (vertex) isoperimetric problem for a graph $G$ is to determine \\ $b_e(i,G)$ ($b_v(i,G)$) respectively for each $i$, $1 \le i \le \vert V \vert$. 
%This problem is NP-hard for a graph.Determining the vertex and edge isoperimetric numbers for even highly structured graphs like hypercubes and grids is difficult. So whenever it is difficult to obtain the actual values of $b_e(i,G)$ and $b_v(i,G)$, we usually  seek upper and lower bounds. Lower bounding vertex isoperimetric numbers is trivial as for all $i$, $b_v(i,G) \geq 1$ and if $mincut(G)$ denotes the value of the mincut for $G$ then for all $i$, $b_e(i,G) \geq mincut(G)$. Note that $b_e(G)$ and $b_v(G)$ act as natural upper bounds. In this paper we study the problem of lower bounding $b_e(G)$ and $b_v(G)$ for trees and for special and interesting cases of trees like complete $t$-ary trees and complete binary trees.
\paragraph{}
Discrete isoperimetric inequalities form a very useful and important 
subject in graph theory and combinatorics. See \cite{Bol86}, Chapter 16 for a  brief introduction on isoperimetric problems. For a detailed treatment see the book by Harper \cite{Har04}. See also the surveys by Leader \cite {Lea91} and by Bezrukov \cite{Bez99,Bez94} for a comprehensive overview of work in the area. The edge(vertex) problem is NP-hard for an arbitrary graph. The NP hardness of the edge version can be seen by observing that if we know $b_e(i,G)$ for all $i$, $1 \leq i \leq \vert V \vert$ we can easily find solutions to the bisection width problem \cite{GJS} and the sparsest cut problem \cite{MatSha}. Isoperimetric problems are typically studied for graphs with special (usually symmetric) structure and the edge and vertex versions of the problem are considered separately as they require different techniques. Probably the earliest example is Harper's work \cite {Har64}: He studied the edge isoperimetric problem for the $d$--dimensional hypercubes. Hart \cite {Hart76} also found the same result separately. Harper later worked on the vertex version \cite {Har66}. Simpler proofs were discovered for his result by Katona \cite {Kat75} and independently by  Frankl and F\"uredi, see \cite {Bol86}, Chapter 16.
The edge isoperimetric problem in the grid i.e. the cartesian product of paths
was considered by Bollabas and Leader \cite {BolLea91}. Since then
many authors have
considered the isoperimetric problems in graph cartesian products. See for example \cite{ChungTetali}.  The isoperimetric problem for the cartesian product of two Markov chains is studied in \cite {HouTetali}. 
Recently Harper considered the isoperimetric problem in Hamming
graphs \cite{Har99}. 
%%%%%Addition
\paragraph{}
The isoperimetric properties of graphs with
respect to  eigen values of their adjacency or Laplacian matrices is
considered by many authors, for example see \cite {AM}. The isoperimetric properties of a graph is very closely related to its expansion properties.
A graph $G$ is called an expander graph if for every 
positive integer $i \le \epsilon |V|$, $b_v(i,G) \ge \epsilon' i$, where
$\epsilon$ and $\epsilon'$ are predefined constants. A great deal of 
effort has gone into explicitly constructing expander graphs - the first construction of an infinite family was due to Margulis \cite{Marg73}. See \cite{ReSaAv} for a recent construction.  
\paragraph*{}
The importance of isoperimetric inequalities lies in the fact that they can be used to give lower bounds for many
useful graph parameters.  
For example it can be shown that $pathwidth(G) \ge b_v(G)$ \cite {ChandranKav05}, $bandwidth(G) \ge b_v(G)$ \cite{Har64} and $cutwidth(G) \ge b_e(G)$ \cite{BezChaHarRS}.  
In \cite {SunCr}, it is shown that  given any $j$ (where  $1 \le j \le \vert V \vert$),
$treewidth(G) \ge   \min_{j/2 \le i \le j} b_v(i,G)-1$ and in \cite{ChandranKav06} it is shown that $carving$-$width(G) \ge \min_{j/2 \le i \le j}  b_v(i,G)$, where $1 \le j \le |V|$ and in \cite{Har64} it is shown that $wirelength(G) \geq \sum_{i=1}^{\vert V \vert}b_e(i,G)$.
\section{Our Results}
Let $T=(V,E,r)$ be a finite,connected rooted tree rooted at $r$. Consider the natural partial order $\preceq_T$ induced by the rooted tree on the vertices.
\begin{definition}
In a rooted tree $T=(V,E,r)$ for any two vertices $u$, $v$, $u \preceq_T v$ if and only if there is a path from the root to $v$ with $u$ in the path. In particular $u \preceq_T u$ for any vertex $u$.
\end{definition}
\begin{definition}
For a rooted tree $T=(V,E,r)$ we define a weight function $w_T:V \rightarrow {1,2,\cdots},\vert V \vert$ as follows: 
$
w_T(u)=\vert \{v \in V:u \preceq v \}\vert 
$(i.e the number of successors of $u$, including $u$). Let us define the weight index of the rooted tree $T=(V,E,r)$ as $\eta(T)=\vert \{w_T(u):u \in V \} \vert$. Note that this is the number of distinct weights. When there is no confusion let $\eta(T)$ be abbreviated by $\eta$.
\end{definition}
\begin{definition}
\label{def2}
For any graph $G$ let, 
$
\ell_G(k) = \vert \{ i \in \mathbb{N}:  1 \leq i \leq \vert V \vert 
\mbox {, }  b_e(i,G) \le k \} \vert
$ where $k$ is a positive integer.
\end{definition}
In other words $\ell_G(k)$ is the number of integers $i$ such that the edge isoperimetric value of $G$ at $i$ is at most $k$. The main Theorem in this paper is as follows: 
%We first give a result on relating $\eta$ and $\ell_T(k)$ for an arbitrary tree.
\begin{theorem}
\label{maintheorem}
  
\end{theorem}
We use the above result to show the following interesting corollaries.
\begin{corollary}
\label{corollary1}
Let $T_d^2$ be the complete binary tree of depth $d$. Then   $c_1d \leq$ \\ $b_e(T_d^2) \leq d$ and $c_2d \leq b_v(T_d^2) \leq d$ where $c_1$ and $c_2$ are constants.
\end{corollary}
%We also get a weaker result for complete $t$-ary trees using similar techniques.
\begin{corollary}
\label{corollary2}
Let $T^t_d$ be the complete $t$-ary tree of depth $d$ with $t \geq 2$ and $d \geq c\log{t}$ where $c$ is a suitable chosen constant. Then,
$c_1\sqrt{t}d \leq b_e(T^t_d) \leq (t-1)d$ and $c_2 \frac{d}{\sqrt{t}} \leq b_v(T^t_d) \leq d$ where $c_1$ and $c_2$ are appropriate constants.
\end{corollary}

%%%%%%
%Addition starts here
%%%%%%
We would like to point that recently Otachi and Yamazaki have considered the problem of determining the vertex isoperimetric peak in complete $t$-ary trees \cite{OatYam}. They prove that $d \geq b_v(T^t_d) \geq  \frac{d\log{t}-(t+6+2\log{d})}{(t+6+2\log{d})}$. Asymptotically our results are better as we prove $b_v(T^t_d) \geq c_2 \frac{d}{\sqrt{t}}$ where $c_2$ is a constant. The best bound that can be obtained from their result for the edge isoperimetric peak is $b_e(T^t_d) \geq  \frac{d\log{t}-(t+6+2\log{d})}{(t+6+2\log{d})}$ while we show that $b_e(T^t_d) \geq c_1\sqrt{t}d $ where $c$ is a constant. Similarly in the special case of a complete binary tree their result implies $b_v(T^2_d) \geq \frac{d\log{2}-(8+2\log{d})}{(8+2\log{d})} \approx \frac{cd}{\log{d}}$ where $c$ is a constant. In contrast we give a tight result showing that $b_e(T_d^2) \geq c_1d$ and $b_v(T_d^2) \geq c_2d$ where $c_1$ and $c_2$ are constants. Morever our proof techniques are such that the above results can be extended to arbitrary (rooted) trees. The proofs in this paper are also comparitively simpler. 
%%%%%
%Addition ends here
%%%%%
As consequences of the above results we have the following theorems. We just mention the theorems here. The necessary definitions and detailed discussions are available in the corresponding sections (section \ref{bssect}-section \ref{essect})
% \begin{theorem}
% \label{theoremwl}
% For a tree of bounded degree $T=(V,E,r)$ with $\vert V \vert=n$ and weight index $\eta \leq c \log{\vert n \vert}$ where $c < \frac{1}{2}$, $wirelength(T)=\Omega(n \log{n})$   
% \end{theorem}
\begin{theorem}
There exists an increasing function $f$ such that for any graph $G$ if $pathwidth(G) \geq k$ then there exists a minor $G'$ of $G$ such that $b_v(G') \geq f(k)$.
\end{theorem}
\begin{theorem}
For the complete binary tree on $T_d^2$ on $n$ vertices $thinness(T_d^2)=\Omega(\log{n})$. This means that there exist trees with arbitrarily large thinness.  
\end{theorem}

\section{Upper bounds on the isoperimetric peak of a tree}

\label{ub}

%\ref{ub} 

% The reader can verify that the folowing assertions are true. For any graph $G$, 
% if $\Delta$ is the maximum degree then, $b_e(i,G) \geq b_v(i,G) \geq \frac{b_e(i,G)}{\Delta}$
% and $b_e(G) \geq b_v(G) \geq \frac{b_e(G)}{\Delta}$. It can also be shown that for any tree $ T=(V,E)$, 
% $pathwidth(T) \leq \log { \vert V\vert}$ \cite{KoSol} and $b_v(T) \leq pathwidth(T)$. Thus we have $b_v(T) \leq \log{\vert V \vert}$ 
% and $b_e(T) \leq \Delta\log{\vert V \vert}$.  So for a complete binary tree of depth $d$ $T_d^2$,
% we have $b_v(T_d^2) \leq d$ and $b_e(T_d^2) \leq cd$ where c is some constant. Similarly for a complete $t$-ary tree of depth $d$ we have,
% $b_v(T_d^t) \leq d)$ and $b_e(T_d^2) = O((t+1)d)$. We can actually show that for $b_e(T_d^t) = O(td)$. Perform an in order 
% traversal of the tree which gives a linear arrangement of the vertices and take the set $S_i$ which consists of the first
% $i$ vertices in the order in which they appear in this arrangement. Clearly $b_e(i,T_d^t) \leq \vert \delta(S_i,T_d^t) \vert$.
% Now it can be easily verfied that $\vert \delta(S_i,T_d^t) \vert=O(td)$ which implies $b_e(i,T_d^t)=O(td)$.

A depth first traversal is one in which all the subtrees of the given rooted tree are recursively visited before visiting the root. Perform such a traversal of the tree and list the vertices in the order in which they appear in the traversal. This gives an ordering of the vertices. Let us choose $S_i$ as the first $i$ vertices as they appear in this ordering. It can be very easily verified that $b_e(i,T) \leq \vert \delta(S_i,G) \vert \leq (\Delta-1) d$ where $d$ is the depth of the tree and $\Delta$ is the maximum degree of a vertex in $T$. Using the same technique we can prove that $b_v(T) \leq d$. For a $t$-ary tree of depth $d$ this implies $b_e(T_d^t) \leq td$ and $b_v(T_d^t ) \leq d$.

\section{Lower bounds on the isoperimetric peak of a tree}

 \begin {definition}
\label {FunctionDef}
 Let  $T=(V,E,r)$ be a rooted  tree with $|V|=n$ and root $r$, 
  and let $S \subseteq V$. 
 Then we define the function $f_{S,T} : E \cup \{ r \} \rightarrow
  \{w_T(u) : u \in V \} \cup \{ 0\}$  as follows: 
  %\begin{linenomath*}
  \begin {eqnarray}
   f_{S,T} (r)  &=& 0  \mbox { {\bf if}  $r \in V - S$} \nonumber \\
                &=& w_T(r) = n \mbox {  {\bf if} $r \in S$ }  \nonumber \\
    f_{S,T}(e) &=& 0  \mbox{ {\bf if}  $e \in E - \delta(S,T)$} \nonumber 
   \end {eqnarray}
   %\end{linenomath*}
 \noindent Finally  if $e = (u,v) \in \delta(S,T)$ ,
   without loss of generality assume that
   $u$ is a child of $v$ in $T$.   
   Then,
  %\begin{linenomath*}
  \begin {eqnarray} 
  f_{S,T}(e) = f_{S,T} (u,v) &=& w_T(u) \mbox { \bf {if}   $u \in S$ } \nonumber \\
                  &=& -w_T(u)  \mbox { \bf {if}  $u \in V-S$} \nonumber 
  \end {eqnarray}
  %\end{linenomath*}
 \end {definition} 

  \begin {lemma}
  \label {FunctionLemma}
    Let $T=(V,E,r)$ be a tree with root $r$ 
     and let $S \subseteq V$. Then,
    $ f_{S,T}(r) + \sum_{e \in E(T) }  f_{S,T}(e)  = |S|$. 
  \end {lemma} 
  \proof { We use induction on the number of vertices $ \vert V \vert = n$.
    For a rooted tree $T'=(V',E',r)$   with  $\vert V' \vert=1$, it 
   is trivial to verify the Lemma.
  Let the Lemma  be true  for any rooted tree $T'=(V'',E'',r'')$ on at most $n-1$
  vertices (where $n \ge 2$)  and for all possible subsets of $V''$.
   Let  $S$ be an arbitrary
   subset of $V$.  Let $v_1,v_2,\cdots,v_k$ be 
  the  children of $r$ in $T$. We denote by $T_i=(V_i,E_i,v_i)$ the subtree
  of $T$ rooted at $v_i$.
  Let $S_i = S \cap V_i$ for $1 \le i \le k$. Also, let
  $f$ denote the function $f_{S,T} : \{r\} \cup E \rightarrow \{w_T(u): u \in V  \} \cup \{0\}$, let 
  $f^i$ denote the function $f_{S_i,T_i}: \{v_i\} \cup E_i
  \rightarrow \{w_{T_i}(u) : u \in V_i \} \cup \{0\} $. By the induction assumption we have, 
% \begin{linenomath*}
   \begin {eqnarray}
    f^i(v_i) +  \sum_{e \in E_i} f^i(e) = |S_i| \mbox { for  } 1 \le i  \le k 
   \end {eqnarray} 
%  \end{linenomath*}
   \noindent Noting that for any edge $e \in E(T) \cap E(T_i)$,
   $f(e) = f^i(e)$  we have:
%   \begin{linenomath*}
    \begin {eqnarray}{}
    \sum_{e \in E(T)} f(e) &=&  \sum_{i=1}^k\sum_{e \in E_i} f^i(e) + \sum_{i=1}^k f(r,v_i)  \nonumber \\ 
                   &=&  \sum_{i=1}^k |S_i| - f^i(v_i)  + \sum_{i=1}^k f(r,v_i)  
     \end {eqnarray}
%    \end{linenomath*}
\noindent By the  definitions of the functions $f$ and $f^i$ (see
  Definition \ref {FunctionDef}) we have: 
%   \begin{linenomath*}	
   \begin {eqnarray}
               f(r,v_i) - f^i(v_i) &=& 0   \mbox {\bf \  if \   $r \in V-S$}   \\
               f(r,v_i) - f^i(v_i) &=& -w_{T_i} (v_i)= -w_T(v_i)  \mbox {\bf \ if \   $r \in S$} 
   \end {eqnarray}
%   \end{linenomath*}
 \noindent Now substituting Equations (3) and (4) in Equation (2),  we get
%   \begin{linenomath*}
   \begin {eqnarray}
    f(r) + \sum_{e \in E} f(e)  &=&  \sum_{i=1}^k |S_i| = |S|  \mbox {\bf \  if $r \in V - S$}  \nonumber
    \end {eqnarray}
%   \end{linenomath*}
    \noindent and 
%  \begin{linenomath*}
    \begin {eqnarray} 
    f(r) + \sum_{e \in E} f(e) &=& \sum |S_i| + w_T(r) - \sum_{i=1}^k w_T(v_i) \nonumber \\
                                  &=& \sum_{i=1}^k |S_i| + 1 =  |S| \nonumber
    \mbox {\bf \  if $r \in S$}
    \end {eqnarray} 
 %  \end{linenomath*}
\noindent as required.
} $\bbox$

We need the following lemma to prove the corollaries of the next theorem.

\begin{lemma}
For any graph $G=(V,E)$, $b_e(G) \geq b_v(G) \geq \frac{b_e(G)}{\Delta}$
\end{lemma}
\begin{proof}
The first part of the inequality is obvious. Let the edge isoperimetric peak occur at $i$ and the vertex isoperimetric peak at $j$. Since $\Delta$ is the maximum degree, $\Delta b_v(i,G) \geq b_e(i,G) = b_e(G) $(Every vertex can have atmost $\Delta$ edges incident on it). But $b_v(G)=b_v(j,G) > b_v(i,G)$. Therefore  $\Delta b_v(G) =\Delta b_v(j,G) \geq b_e(G)$.
\end{proof}

%\begin{theorem}
%\label{theorem1}
%Let $T=(V,E,r)$ be a rooted tree with $\vert V \vert = n$ and weight %index $\eta$ and let $T=(V^{+},E^{+},\hat{v})$ be an extension of it . %Then
%$
%\ell(k) \leq 2{{2\eta+k} \choose k}
%$.
%\end{theorem}
\noindent  {\bf Theorem \ref {maintheorem}}. {\it For any rooted tree $T=(V,E,r)$, with weight index $\eta$ , $\ell_T(k) \leq 2{{2\eta+k} \choose k}$} 
\begin{proof}
Let $i \le \vert V \vert$ be a positive  integer  such that
  $b_e(i, T)= k' \le  k$.  Then there exists a subset 
  $S_i \subseteq V$ such that $|\delta(S_i,T)| = k'$ and $\vert S_i \vert = i$. 
  Let $\delta(S_i,T) = \{ e_1, e_2, \cdots, e_{k'} \}$. We define
  $k+1$ variables $t_0,t_1, \cdots, t_k$ as follows.  
  Let $t_0 = f_{S_i,T}(r)$ and  let $t_i = f_{S_i,T}(e_i)$ for
  $1 \le i \le k'$.
  If $k' < k$, then let $t_i = 0$ for $k' < i \le k$. By Lemma  
  \ref {FunctionLemma}, we have 
  $\sum_{e \in E} f_{S_i,T}(e) = |S_i| = i$. 
  Recalling Definition \ref{FunctionDef}, for an edge $e$, 
   $f_{S_i,T}(e) \ne 0$
  only when $e \in \delta(S_i, T)$. Thus we have:
%  \begin{linenomath*}
   \begin {eqnarray}
   \label {SumOfTermsEquation}
   t_0 + t_1 + \cdots + t_k = i  \nonumber
   \end {eqnarray} 
%  \end{linenomath*}
  
  How many distinct positive integers can be expressed as $\sum_{i=0}^k t_i$? 
  This will clearly give an upper bound for $\ell(k)$. 
 Let $W=\{w_1,\cdots,w_\eta\}$ where $\eta$ is the weight index of the tree, denote the set of distinct weights. Then $t_i$ can take the values $0$ or $\pm w_j,1 \leq j \leq \eta$. Considering  the $k$ variables $t_i$ ($1 \le i \le k$)
    as $k$ unlabeled
  balls  and imagining the $2\eta+1$ distinct possible values they    
   can take as  $2\eta+1$ labeled boxes, it is easy to see that
   the number of distinct integers expressible as $\sum_{i=1}^k t_i$ is
    bounded  above by the number of ways of arranging $k$ unlabeled balls 
 in $2\eta+1$ labeled boxes, i.e. ${{2\eta+ k} \choose k }$.  
  Recalling that $t_0$ can take only two possible values, we get:  
%\begin{linenomath*}
   $$ \ell_T(k) \le   2 { {2\eta+k} \choose k } $$
%\end{linenomath*}
 \bbox
\end{proof}

\begin{definition}
\label{def12}
For any graph $G$ with weight index $\eta$, define $p$ as the minimum value of k such that $2{{2\eta+k} \choose k} \geq n$
\end{definition}
\begin{lemma}
For any rooted tree $T=(V,E,r)$, $b_e(G) \geq p$
\end{lemma}
\begin{proof}
Assume $b_e(G) < p$. Let $b_e(G)=q$. Then by the definition of $p$ we have $2{{2\eta+q} \choose q} < n$. But by Definition \ref{def2} $\ell_T(q)=n$ a contradiction.
\end{proof} 
\noindent  {\bf Corollary \ref {corollary1}}. { \it Let $T_d^2$  be the complete binary tree of depth $d$. Then   $c_1d \leq$ \\ $b_e(T_d^2) \leq d$ and $c_2d \leq b_v(T_d^2) \leq d$ where $c_1$ and $c_2$ are constants.}
\begin{proof}
Let the number of vertices in $T_d^2$ be denoted by $n$. We need only prove that $b_e(T^2_d) \geq c_1d$ for some constant $c_1$ as the upper bound follows from Section \ref{ub}. Note that $\eta(T^2_d) =d$ so, $\ell(k) \leq 2{{2d+k} \choose k}$ where $k$ is a positive integer. Now let $k = \lfloor \frac{d}{5} \rfloor = \lfloor 0.2d \rfloor$.  Then we have
  (discarding the floor symbol),
%\begin{linenomath*}
\begin{eqnarray}
    2{ {2d + k} \choose k } & = & 2{{2.2d} \choose {0.2d}} \nonumber \\
    & = &  \frac {2(2.2d) !}{0.2d! 2d!} \nonumber \\ 
    & = & \frac {c}{\sqrt d} \left ( \frac {(2.2)^{2.2}} {(0.2)^{0.2} 2^2} \right )^d \nonumber \\
    & \le & \frac {c'}{\sqrt d} (1.96)^d  \nonumber
\end{eqnarray}
%\end{linenomath*}
Here we have used Stirling's approximation,   $c''\sqrt{2\pi n}n^ne^{-n} \leq n! \leq \\ c'''\sqrt{2 \pi n} n^n e^{-n}$. This means that for a sufficiently large value of $d$, $\ell(k) < n$ when $k=\frac{d}{5}$ which implies that $b_e(T^2_d) \geq c_1d$ where $c_1$ is a constant. Again this implies $b_v(T^2_d) \geq c_2d$ where $c_2$ is a constant, as $\Delta=3$ for a complete binary tree.   \bbox  
\end{proof}
The reader may note that the above proof shows that for almost all integers $i$, $1 \leq i \leq k$ $b_e(i,T_d^2) \geq .2d$. More precisely $\lim_{d \rightarrow \infty} \frac{\ell_{T_d^2}(\frac{d}{5})}{n} \rightarrow 0$. \\
%We can extend the above result to complete $t$-ary trees.\\ 
%\begin{corollary}
%Let $T^t_d$ is a complete $t$-ary tree of depth $d$. Then,
%$b_e(T^t_d) = \omega(\sqrt{t}d)$, %$b_v(T^t_d)=\omega(\frac{d}{\sqrt{t}}$ and $b_e(T^t_d)=O(td)$, %$b_v(T^t_d)=O(d)$
%\end{corollary}
\noindent  {\bf Corollary \ref {corollary2}}. {\it Let $T^t_d$ be the complete $t$-ary tree of depth $d$ with $t \geq 2$ and $d \geq c\log{t}$ where $c$ is a suitable chosen constant. Then,
$c_1\sqrt{t}d \leq b_e(T^t_d) \leq td$ and $c_2d \frac{d}{\sqrt{t}} \leq b_v(T^t_d) \leq d$ where $c_1$ and $c_2$ are constants.}
\begin{proof}
The upper bound follows from Section \ref{ub}. We will assume that $t \geq 9$ initially and $d \geq 30$. Note that for a $t$-ary tree of depth $d$, $\eta(T_d^t) = d$. For a positive integer $k$, by Theorem \ref{maintheorem} we have 
%\begin{linenomath*}
$$\ell(k) \leq 2{{2d+k} \choose k}$$
%\end{linenomath*}
Now let $k = \lfloor m\sqrt{t}d \rfloor$ where $0 < m < 2(\frac{1}{e}-\frac{1}{3})$ is a constant. Then we have(discarding the floor symbol),
 %\begin{linenomath*}
\begin{eqnarray}
    2{ {2d + k} \choose k } & = & 2{{(2+m\sqrt{t})d} \choose {m\sqrt{t}d}} \nonumber \\
    & = &  \frac {2((2+m\sqrt{t})d)!}{(m\sqrt{t}d)! (2d)!} \nonumber \\ 
    & \leq & \frac {c'\sqrt{(2+m\sqrt{t})d}}{\sqrt{2d}\sqrt{m\sqrt{t}d}} 
  \left ( \frac {(2+m\sqrt{t})^{2+m\sqrt{t}}} {(m\sqrt{t})^{m\sqrt{t}} 2^2} \right )^d \nonumber \\
    & \leq & c'' \left ( \frac {(2+m\sqrt{t})^{2+m\sqrt{t}}} {(m\sqrt{t})^{m\sqrt{t}} 2^2} \right )^d
\end{eqnarray}
%\end{linenomath*}
as $\frac {\sqrt{(2+m\sqrt{t})d}}{\sqrt{2d}\sqrt{m\sqrt{t}d}} < 1$ for $d \geq 30$ and $t \geq 9$ with $m$ being chosen appropriately. Now consider,  
%\begin{linenomath*}
\begin{eqnarray}
c'' \left ( \frac {(2+m\sqrt{t})^{2+m\sqrt{t}}} {(m\sqrt{t})^{m\sqrt{t}} 2^2} \right )^d 
& = & c'' \left ( \frac{(2+m\sqrt{t})^{m\sqrt{t}}}{(m\sqrt{t})^{m\sqrt{t}}}\frac{(2+m\sqrt{t})^2}{2^2} \right )^d \nonumber \\
& = & c'' \left (( (1+\frac{2}{m\sqrt{t}})^{\frac{m\sqrt{t}}{2}})^{2}\frac{(2+m\sqrt{t})^2}{2^2} \right )^d \nonumber \\
& \leq & c'' \left (e^2 \frac{(2+m\sqrt{t})^2}{2^2} \right )^d 
\end{eqnarray}
%\end{linenomath*}
Here we have used the fact that $(1+x)^{\frac{1}{x}} \leq e$ for $x > 0$. Let the number of nodes in $T_d^t$ be $n = \frac{(t^d-1)}{(t-1)} \geq t^{(d-1)}$. Therefore from Eqns (5) and (6) we have
%\begin{linenomath*}
\begin{eqnarray}
\frac{2{ {2d + k} \choose k }}{n}  
& \leq & \frac{2{ {2d + k} \choose k }}{t^{(d-1)}}  \nonumber \\
& \leq & \frac{c'' \left ( e^2 \frac{(2+m\sqrt{t})^2}{2^2} \right ) ^d}{t^{(d-1)}} \nonumber \\
& \leq & c'' t \left ( e^2 (\frac{1}{\sqrt{t}}+\frac{m}{2})^2 \right ) ^d \nonumber \\
& = & S(say)  \nonumber 
\end{eqnarray}
%\end{linenomath*}
Clearly for large enough $d$ i.e $d \geq c\log{t}$, $S < 1$ as $ e^2 (\frac{1}{\sqrt{t}}+\frac{m}{2})^2 < 1$ for the chosen value of $m$ . This means that for large enough $d$, ${ {2d + k} \choose k } < n$ which implies $b_e(T_d^t) \geq p \geq k \geq m\sqrt{t}d$. In our proof we have assumed that $t \geq 9$. This assumption can be removed by noting that for all values of $t < 9$ we can prove $b_e(T_d^t) > c'''d$ for some constant $c'''$ using the same techniques as in the proof for the binary tree . So we can show $b_e(T_d^t) \geq c''''\sqrt{t}d$ by taking $c'''' = \frac{c'''}{3}$ since in this case $\sqrt{t} < 3$. This completes the proof that $b_e(T_d^t) \geq c_1\sqrt{t}d$ for all $t \geq 2$ where $c_1$ is a appropriately chosen constant. $\Delta=(t+1)$ in $T_d^t$. Therefore $b_v(T_d^t) \geq \frac{b_v(T_d^t)}{(t+1)} \geq  \frac{c_1\sqrt{t}d}{(t+1)} \geq \frac{c_2d}{\sqrt{t}}$. \bbox  
%  It is easy to check that $\lim_{d \rightarrow \infty} S \rightarrow 0$ if $m < 2(\frac{1}{e}-\frac{1}{3})$. Choosing this value of $m$ and Eqns (5), (6) and (7) lead us to conclude that $k  \geq c'\sqrt{t}d$ which in turn implies that $b_e(T_d^t) = c_1\sqrt{t}d$ and $b_v(T_d^t) = c_2\frac{d}{\sqrt{t}}$ where $c'$, $c_1$ and $c_2$ are constants 
\end{proof}
These results can be generalized to an arbitrary tree.

\begin{corollary}
 Let $T=(V,E,r)$ be a rooted tree with $\vert V \vert =n$ and weight index $\eta$ and $p \geq 2$. Then, $b_e(T) \geq c_1\eta(n^{(\frac{1}{2\eta})}-c_2)$ and $b_v(T) \geq \frac{c_1\eta(n^{(\frac{1}{2\eta})}-c_2)}{\Delta}$ where $c_1$ and $c_2$ are constants.
\end{corollary}
\begin{proof}
%\begin{proof}
%We have from Definition \ref{def12}
%$$
We have, $n \leq {{2\eta+ p} \choose p} $. 
%$$
Let $p=\omega\eta$. Then, 
%\begin{linenomath*}
\begin{eqnarray}
n & \leq & 2 {{2\eta+\omega\eta } \choose \omega\eta} \nonumber \\
  & = & 2 {{(2+\omega)\eta } \choose \omega\eta} \nonumber
\end{eqnarray}
%\end{linenomath*}
%Now using Sterling's approximation 
%$
%c_{1}\sqrt{2\pi n}n^{n}e^{-n} \leq n!\leq c_{2}\sqrt{2\pi n}n^{n}e^{-n}
%$
%where $c_{1}$ and $c_{2}$ are two suitable constants.,
%\begin{linenomath*}
\begin{eqnarray}
n & \leq & \frac{2c\sqrt{2\pi(2+\omega)\eta)}((2+\omega)\eta)^{(2+\omega)\eta}e^{-(2+\omega)\eta}}{(c'\sqrt{2\pi\omega\eta}(\omega\eta)^{\omega\eta}e^{-\omega\eta})(c'\sqrt{4\pi\eta}(2\eta)^{2\eta}e^{-2\eta})} \nonumber \\
& \leq & \frac{c''((2+\omega)\eta)^{(2+\omega)\eta}}{((\omega\eta)^{\omega\eta})((2\eta)^{2\eta})} \nonumber \\ & \leq & c'''(1+\frac{2}{\omega})^{\omega\eta}(\frac{\omega}{2}+1)^{2\eta} \nonumber
\end{eqnarray}
%\end{linenomath*}

where $c$, $c'$, $c''$ and $c'''$ are suitably chosen constants. We have used the fact that $2\eta+\omega\eta \leq 2\eta\omega\eta$ (which follows from the fact that $2\eta \geq 2$ and $p = \omega\eta \geq 2$).  Simplifying this yields,
%\begin{linenomath*}
$$
n^{\frac{1}{2\eta}} \leq c''''(1+\frac{2}{\omega})^{\frac{\omega}{2}}(\frac{\omega}{2}+1)
$$
%\end{linenomath*}
Since $(1+x)^{\frac{1}{x}} \leq e$ for $x > 0$,
%\begin{linenomath*}
$$
c_1n^{\frac{1}{2\eta}}-2 \leq \omega
$$
%\end{linenomath*}
Since $b_e(T) \geq p$, $b_e(T) \geq c_1\eta(n^{(\frac{1}{2\eta})}-c_2)$ where $c_1$ and $c_2$ are suitably chosen constants. Therefore $b_v(T) \geq \frac{c_1\eta(n^{(\frac{1}{2\eta})}-c_2)}{\Delta}$ and the corollary follows. \bbox
%\end{proof}
\end{proof}

\noindent {\bf Comment: }
It is interesting to study for what values of $\eta$ the above result would be useful. A simple obsevation is that $n^{(\frac{1}{2\eta})}>c_2$. An analysis of the proof for the above result shows that that $c_1 \geq \frac{1}{e}$ and thus $c_2 \leq 2e$. We note that $n=e^{\log{n}}$. For a tree $T$ with $ \eta \leq \frac{\log{n}}{4}$ we would have $b_e(T) \geq c\eta$ and $b_v(T) \geq \frac{c\eta}{\Delta}$ for a constant $c$. Similarly for a tree $T$ with $\eta=k$ a constant we have $b_e(T) \geq c'n^{\frac{1}{2k}}$ and $b_v(T) \geq \frac{c'n^{\frac{1}{2k}}}{\Delta}$ for a constant $c'$.

\section{Applications}

\subsection{Pathwidth} 
\label{bssect}
Pathwidth and Path decomposition are important concepts in graph theory and computer science. For the definition and several applications see \cite{Bodland3}. It is not difficult to show that $pathwidth(G) \geq b_v(G)$ 
(see \cite{ChandranKav05}). An obvious question is whether the reason for
 the high pathwidth of a graph $G$,  is the  ``good'' isoperimetric property
 of an   induced subgraph or minor of $G$ . More precisely if 
$pathwidth(G) \geq k$ is it possible to find an induced subgraph 
or minor $G'$ of $G$ such that $b_v(G') \geq f(k)$ for some 
function $f$, where $f(k)$ increases with $k$. 
 Let us first consider whether such an  induced subgraph
 always exists.  The answer is negative: Given any integer $k$,
 it is possible to demonstatrate
 a graph $G$ (on arbitrarily large number of vertices) such that 
 pathwidth$(G) \ge k$,  but  $b_v(G')$ for any induced subgraph of $G$
 is  bounded above by a constant. For example, one can start with
 a complete binary tree of sufficiently large depth. The pathwidth of
 such a tree is $\Omega(d)$, where $d$ is the depth. Now we can replace
 each edge of the binary tree with a path of appropriately chosen length,
 to make sure that for any induced subgraph $T'$  of the resulting tree
 $b_v(T')  \le c$, where $c$ is some constant. On the other hand,
 reader can easily verify that  by replacing an edge with 
  a path (i.e. by subdividing an edge)  we can not decrease the pathwidth
 of the original graph. Thus the resulting tree will have pathwidth as
 much as that of the original. (We leave the rigorous proof of the
 above as an exercise to the reader.) 
But when we ask the same question  with respect to 
 minors,  the answer is positive.
Robertson and Seymour proved the following result.
(See \cite {diest}, Chapter 12.)
 {\it If pathwidth$(G) \ge k$, then there
exists a function $g$ such that every tree on at most $g(k)$ vertices is a
minor of $G$. }
Then clearly there exists a minor of $G$ which is isomorphic to a  
complete binary tree 
$T$ on at least $\frac {g(k)}{2}$ vertices. 
By our  result (Corollary \ref {corollary1}) $b_v(T) \ge  c\log {n} =
c'  \log {(g(k))}$  where $c$ and $c'$ are appropriate  constants.
Thus we have the following result: 

%\begin {theorem}
\textit{There exists a function $f$  such that if the pathwidth of a graph $G$
 is at least $k$, then there exists a minor $G'$ of $G$ such that
$b_v(G') \ge f(k)$.}

\subsection{Thinness} 
\label{essect}
A new graph parameter {\it thinness}, is defined in \cite {ManOrRiCha}
 which attempts to generalize certain properties of interval graphs.
  The thinness of a graph $G=(V,E)$ is the minimum
positive integer $k$ such that there exists an ordering 
$v_1,v_2, \cdots, v_n$ (where $n = |V|$) of the  vertices of $G$  
  and a partition $V_1,V_2, \cdots, V_k$ of $V$ into $k$ disjoint
sets, satisfying the following condition:
For any triple $(r,s,t)$ where $r < s < t$, if $v_r$ and $v_s$ 
belong to the same set $V_i$ and if $v_t$  is adjacent to $v_r$
then $v_t$ is adjacent to $v_s$ also.  The motivation for studying 
this parameter was the observation  that the maximum independent set 
problem  can be solved in polynomial time, if a family of graphs has  bounded
thinness. The applications of {\it thinness} for the Frequency Assignment 
Problems in GSM networks are explained in \cite {ManOrRiCha}. One intersting aspect of thinness is that for a graph $G$, $thinness(G) \leq pathwidth(G)$. 
A natural question which arose in connection with
our study of thinness was  the following: {\it Are trees of bounded thinness?
In other words, is there a family of trees for which the
 thinness grows with the number of vertices?} It is 
proved in a later paper by the authors of \cite{ManOrRiCha}  that for any graph  $G$, 
thinness$(G) \ge  \frac {b_v(G)}{\Delta}$ where
$\Delta$ is the maximum degree of $G$.  Combining this lower bound
with our earlier observations, 
 we can infer that the thinness of  a complete binary tree on $n$ vertices 
is $\Omega (\log n)$.

\end{document}